\documentstyle[amssymb,11pt]{article}

\newtheorem{defi}{Definition}[section]
\newtheorem{theo}{Theorem}[section]
\newtheorem{prop}{Proposition}[section]
\newtheorem{lemma}{Lemma}[section]
\newtheorem{rem}{Remark}[section]
\newtheorem{cor}{Corollary}[section]

\newcommand{\bdef}{\begin{defi}}
\newcommand{\ede}{\end{defi}}
\newcommand{\bsat}{\begin{theo}}
\newcommand{\esat}{\end{theo}}
\newcommand{\bprop}{\begin{prop}}
\newcommand{\eprop}{\end{prop}}
\newcommand{\blem}{\begin{lemma}}
\newcommand{\elem}{\end{lemma}}
\newcommand{\brem}{\begin{rem}}
\newcommand{\erem}{\end{rem}}
\newcommand{\bcor}{\begin{cor}}
\newcommand{\ecor}{\end{cor}}
\newcommand{\bbew}{{\bf Proof.}\hspace{0.5cm}}
\newcommand{\ebew}{\hspace{0.5cm}\bigskip\mbox{$\Box$}}

\newcommand{\be}{\begin{equation}}
\newcommand{\ee}{\end{equation}}
\newcommand{\bib}{\bibitem}

\newcommand{\hra}{\hookrightarrow}
\newcommand{\ra}{\rightarrow}
\newcommand{\ts}{\textstyle}
\newcommand{\A}{{\cal A}}

\newcommand{\N}{{\Bbb N}}
\newcommand{\R}{{\Bbb R}}
\newcommand{\C}{{\Bbb C}}
\newcommand{\Z}{{\Bbb Z}}

\newcommand{\K}{{\frak k}}
\newcommand{\G}{{\frak g}}
\newcommand{\T}{{\frak t}}

\newcommand{\V}{{\cal V}}
\newcommand{\D}{{\d_\G}}
\renewcommand{\L}{{\cal L}}
\renewcommand{\O}{{\cal O}}

\renewcommand{\d}{{\mbox{d}}}

\newcommand{\f}{\frac}
\newcommand{\F}{{\cal F}}
\newcommand{\p}{\partial}
\renewcommand{\l}{\langle}
\newcommand{\r}{\rangle}

\title{On Riemann-Roch Formulas for Multiplicities}
\author{\bigskip  Eckhard Meinrenken \thanks{
Mathematics Department, M.I.T., Cambridge MA 02139, U.S.A.
Address after August 1, 1994:
Fakult\"at f\"ur Physik, Hermann-Herder-Str. 3,
D-79104 Freiburg, Germany.}\,\,\thanks{
Supported by a grant from the
German Academic Exchange Service. }}

\date{May 1994 \thanks{Revised Version, October 1994}}

\begin{document}
\maketitle
\begin{abstract}
A Theorem due to Guillemin and Sternberg \cite{GS82a}
about geometric quantization
of Hamiltonian actions of compact Lie groups $G$
on compact K\"ahler manifolds
says that the dimension of the $G$-invariant subspace is
equal to the Riemann-Roch number of the  symplectically reduced space.
Combined with the shifting-trick, this gives explicit formulas
for the multiplicities of the various irreducible components.
One of the assumptions of the Theorem is that the reduction
is regular, so that the reduced space is  a smooth symplectic manifold.
In this paper,
we prove a generalization of this result to the case where the reduced
space may have orbifold singularities.
Our proof uses
localization techniques from equivariant cohomology, and relies in
particular on recent work of Jeffrey-Kirwan \cite{JK93} and
Guillemin \cite{G94}. Since there are no complex geometry arguments
involved, the result also extends to non K\"ahlerian settings.
\end{abstract}

\section{Introduction}

Let $(M,\omega)$ be a compact K\"{a}hler manifold, and
let $\tau:L\ra M$ be a holomorphic line bundle over $M$
with Hermitian fiber metric $h$. $(L,h)$ is said to satisfy
the quantizing condition if $-2\pi i \omega$ is the curvature of
the canonical Hermitian connection on $L$.
Let $H^i(M,{\cal O}(L))$ be the $i$th cohomology group for the
sheaf of germs of holomorphic sections. By the
Riemann-Roch Formula of Hirzebruch and Atiyah-Singer\cite{AS68a},
the Euler number
\be \mbox{Eul}(L):=\sum_i (-1)^i \dim H^i(M,{\cal O}(L))\ee
is equal to the characteristic number
\be \mbox{Eul}(L)=\int_M Td\,(M)\,Ch\,(L),\label{HRR}\ee
where ${Td}\,(M)$ is the Todd class and ${Ch}\,(L)=e^{[\omega]}$
the Chern character.
Recall that if $L$ is ``sufficiently positive'', in particular
if one replaces $L$ by some sufficiently high tensor power,
all the cohomology groups with $i>0$ are zero by
Kodaira's Theorem \cite{GH78}, so in this case
(\ref{HRR}) gives a formula for the dimension of the space
$H^0(M,{\cal O}(L))=\Gamma_{hol}(M,L)$ of holomorphic sections.

Let $G$ be a compact, connected Lie group that acts on $M$
by K\"ahler diffeomorphisms $\Phi:G\times M\ra M$,
with an equivariant moment map $J:M\ra \G^*$.
Suppose also that $\Phi$ lifts to Hermitian bundle automorphisms
of $L\ra M$, according to the rules of geometric
quantization \cite{GS82a}.
The corresponding virtual representation of $G$ on
$\sum (-1)^i\,H^i(M,{\cal O}(L))$ may then be regarded as the
``quantization'' of the classical action $\Phi$.
Its character $\chi$ is the element of the representation ring
$R(G)$ defined by
\be \chi(g):=\sum_i\, (-1)^i\, \mbox{tr}\big(g|H^i(M,{\cal O}(L))\big).
\label{char}
\ee
{}From the Equivariant Riemann-Roch Formula of Atiyah-Segal-Singer
\cite{AS68b,AS68a}, one has an expression for $\chi(g)$
as the evaluation of certain characteristic classes on
the fixed point set $M^g=\{x\in M|\,g.x=x\}$ (which is a
K\"ahler submanifold of $M$):
Let ${Ch}^g(L|M^g)=c_L(g)\,{Ch}(L|M^g)$, where $c_L(g)\in S^1$ is the
(locally constant) action of $g$ on $L|M^g$. Denote by $N^g$ the
normal bundle of $M^g$ in $M$, by $F(N^g)$ its curvature, and let
\be D^g(N^g)=\det(I-(g^\sharp)^{-1}e^{-\f{i}{2\pi}F(N^g)})\ee
where $g^\sharp$ is the automorphism of $N^g$ defined by $g$.
Then
\be \chi(g)=\int_{M^g} \frac{{Td}\,(M^g){Ch}^g(L|{M^g})}{D^g(N^g)}.
\label{character}\ee

By a Theorem of Guillemin and Sternberg \cite{GS82a}, there are also
Riemann-Roch Formulas for the multiplicities of the
irreducible components of the above representation,
at least if certain regularity
assumptions are satisfied. Let $T\subset G$ be a maximal torus, and
$\G=\T\oplus [\T,\G]$ the corresponding decomposition of the Lie
algebra.
Choose a fundamental Weyl chamber
$\T^*_+\subset\T^*\subset \G^*$, let $\Lambda\subset \T^*$ the
integral lattice, and $\Lambda_+=\Lambda\cap \T^*_+$ the dominant
weights.
For a given lattice point $\mu\in\Lambda_+$,
let $V_\mu$ denote the corresponding
irreducible representation with highest weight $\mu$,
and define the multiplicity
$N(\mu)$ by the alternating sum
\be N(\mu):=\sum_i (-1)^i \dim\big(V_\mu^*\otimes
H^i(M,{\cal O}(L))\big)^G. \label{defmult}\ee

Suppose that $\mu\in\Lambda_+$ is a regular value of $J$, or equivalently
that the action of  the isotropy group
$G_\mu$ on $J^{-1}(\mu)$ is locally free.
If the action is in fact free, the reduced space
$M_\mu=J^{-1}(\mu)/G_\mu$ is a smooth symplectic manifold,
and it is well-known that it acquires a natural K\"ahler
structure, together with a quantizing line bundle $L_\mu$.
The main result of \cite{GS82a} is that the multiplicity of $\mu$ in
$\Gamma_{hol}(M,L)$ is equal to the dimension
of the space $\Gamma_{hol}(M_\mu,L_\mu)$, so in particular $N(\mu)$ is
given by the Euler number of $L_\mu$ if $L$ is sufficiently positive:

\bsat[V. Guillemin, S. Sternberg \cite{GS82a}] If the action of
$G_\mu$ on $J^{-1}(\mu)$ is free, and if $L$ is sufficiently
positive,
\be N(\mu)=\int_{M_\mu}{Td}\,(M_\mu){Ch}\,(L_\mu).  \label{GS}    \ee
\esat
The ``physical'' interpretation of this Theorem is that reduction and
quantization commute.

In practice, one is often dealing with situations where the action is
only locally free. In this case, the reduced
space is in general just an
orbifold (or V-manifold) in the sense of Satake \cite{S57},
which means (roughly) that it is locally the quotient of a manifold
by a finite group.
Moreover, the reduced line bundle
is in general just an orbifold-bundle, that is, at some points the
fiber of $L_\mu$ is not $\C$, but its quotient by a
finite group.
Guillemin and Sternberg conjectured that in this case,
the right hand side of
(\ref{GS}) has to be replaced by the expression from T. Kawasaki's
Riemann-Roch Formula for orbifolds  \cite{K79}.
It was proved by R. Sjamaar \cite{S93} that this assertion is true if
$L$ is sufficiently positive, and if $L_\mu$ is an honest line bundle.
In fact, his approach also covers the truly singular case where
$\mu$ is not even a regular value, by using Kirwan's partial
desingularization procedure to reduce it to the orbifold case.
On the other hand, the condition that $L_\mu$ be a genuine line
bundle
is rather restrictive. It is the aim of the present paper
to give a different proof of the Guillemin-Sternberg conjecture (for
$\mu$ a regular value), without having to make this assumption.

The method we use is motivated by recent work of V. Guillemin
\cite{G94},
who used localization techniques from equivariant cohomology
to establish the connection between the Multiplicity Formula (\ref{GS})
and a certain formula for counting lattice points in polytopes.
This formula is known to be true in various interesting cases, and
for these gives a new proof of (\ref{GS}) without using any
complex geometry arguments. (In particular, it also works for
{\em almost} K\"ahler polarizations.)

We will adopt this utilization of equivariant cohomology,
but in a slightly different guise.
The main idea is to consider the rescaled problem, where we replace
\be \omega\leadsto m\omega,\,L\leadsto\,L^m,\,J\leadsto\, mJ,
\,\mu\leadsto m\mu\ee
for $m\in\N$.

Our starting point will be the Equivariant Riemann-Roch Formula,
but instead of (\ref{character}) we will use it in a
form due to  Berline and Vergne \cite {BV85},
involving equivariant characteristic classes.
By a stationary phase
version of the Localization Formula of Jeffrey and Kirwan \cite{JK93},
we pass from equivariant characteristic classes to (ordinary)
characteristic classes on the reduced spaces. This leads to
the desired Multiplicity Formula for $N^{(m)}(m\mu)$,
up to an error term $O(m^{-\infty})$.

Since the multiplicities are integers, one easily finds
that the error term is zero for large $m$.
To investigate the general dependence of $N^{(m)}(m\mu)$ on $m$, we
use a different expression for $N^{(m)}(m\mu)$ via
the number of lattice points in certain polytopes.
If $J(M)$ is contained in the set of regular elements of $\G^*$,
in particular in the abelian case, this analysis
turns out to be sufficiently
good to show that the above error term is zero for all $m$.

\section{Statement of the result}

In order to state the result, we have to give a closer
description of the reduced space and its singular strata.
Suppose that $\mu\in\Lambda_+$ is a regular value of $J$.
Recall first the shifting trick to express
$M_\mu$ as a reduced space at the zero level set:
Let ${O}=G.\mu$ be the coadjoint orbit through $\mu$,
equipped with its usual Kirillov K\"ahler structure, and let
${O}^-$ denote ${O}$ with the opposite K\"ahler
structure. The action of $G$ on ${O}$ is Hamiltonian, with
moment map $\Psi$ the embedding into $\G^*$.
Then $\tilde{M}=M\times {O}^-$ is a K\"ahler manifold,
and the diagonal
action of $G$ is Hamiltonian, with moment map $\tilde{J}=J-\Psi$.
There are canonical identifications
\be M_\mu=J^{-1}(\mu)/G_\mu\cong J^{-1}({O})/G\cong
   \tilde{J}^{-1}(0)/G.\ee
By Kostant's version of the Borel-Weil-Bott Theorem,
one also has a natural quantizing bundle $\Xi\ra {O}$, and
the irreducible representation $V_\mu$ corresponding to $\mu$
gets realized as the space of holomophic sections
of $\Xi$. (The higher order cohomology groups $H^i(O,\O(\Xi))$
vanish.) The tensor product
$\tilde{L}:=L\otimes \Xi^*$ quantizes
$\tilde{M}$, and there is an isomorphism
\be H^i(\tilde{M},{\cal O}(\tilde{L}))\cong
V_\mu^*\otimes H^i(M,{\cal O}(L)).\label{qst}\ee
Hence $N(\mu)={\rm Eul}(\tilde{L})$, which is the quantum counterpart
of the shifting-trick.

Using the shifting-trick, it is enough to consider the case $\mu=0$.
The reduced space $M_0$ inherits a natural K\"ahler structure
from $M$ (see \cite{GS82a}), and the reduced bundle $L_0=(L|J^{-1}(0))/G$
renders a quantizing orbifold-line bundle. Note however
that $L_0$ need not be an honest line bundle, not even
over the smooth part of $M_0$.
Sections of an orbifold bundle are defined as coming from
invariant sections for the local orbifold charts,
so all sections of $L_0$ have to vanish at points were the fiber
is not $\C$.

Let us regard $P=J^{-1}(0)$ as an orbifold-principal bundle over
$M_0=J^{-1}(0)/G$. Following \cite{K79,F92}, we introduce
\be \tilde{P}=\{(x,g)|\,x\in P,\,g.x=x\}\subset P\times G,\ee
and let $\Sigma=\tilde{P}/G$ be its quotient under the locally free action
$h.(x,g)=(h.x,\,h\,g\,h^{-1})$.
The projection of $\tilde{P}$ to the second factor descends to
a locally constant mapping
\be \tau:\Sigma\ra {\rm Conj}(G)\ee
to the set of conjugacy classes. For $g\in G$, let $(g)={\rm Ad}(G).g$
denote the corresponding conjugacy class, and $\Sigma_g$ its preimage under
$\tau$. There is a natural identification $\Sigma_g=P^g/Z_g$, where
$P^g\subset P$ is the fixed point
manifold and $Z_g$ the centralizer of $g$ in $G$.
Since the fixed point set $M^g\subset M$ is a K\"ahler submanifold, and
the action of $Z_g$ on $M^g$ is Hamiltonian with the restriction of $J$
serving as a moment map, this makes clear that $\Sigma$ is a
K\"ahler orbifold (with several components of different dimensions).
Note that this K\"ahler structure does not depend on the choice of the
representative for $(g)$. Observe also that $\Sigma_e\cong M_0$.

The collection of bundles $(L|P^g)/Z_g$ defines a
quantizing orbifold line bundle $L_\Sigma\ra \Sigma$.
As above, let $c_L(g)\in S^1$ be the locally constant action of
$g$ on $L|P^g$, denote by $c_\Sigma:\Sigma\ra S^1$ the
function defined by the $c_L(g)'s$, and let ${Ch}^\Sigma(L_\Sigma)$
be the cohomology class defined by
\be Ch^\Sigma(L_\Sigma)=c_\Sigma\,e^{\omega_\Sigma}\ee
where $\omega_\Sigma$ is the K\"ahler form on $\Sigma$.

Consider now the natural mapping $f:\Sigma\ra M_0$, sending
$G.(x,g)\ra G.x$. In a local orbifold chart,
the tangent space to $\Sigma$ at $G.(x,g)$ is isomorphic to
$T_x(M^g\cap J^{-1}(0))/T_x(Z_g.x)$, while the tangent space to
$M_0$ at $G.x$ is $T_x(J^{-1}(0))/T_x(G.x)$. From this,
it is easy to see that $f$ is a Khler immersion.
Let $N_\Sigma\ra \Sigma$ be the normal bundle of this immersion,
and denote by $g^\sharp$ the automorphism of $N_\Sigma|\Sigma_g$
induced by the action of $g$.
Then the collection of differential forms
\be \det(I-(g^\sharp)^{-1}e^{-\f{i}{2\pi}F(N_\Sigma)}),\ee
where $F(N_\Sigma)$ is the curvature of $N_\Sigma$, defines a characteristic
class $D^\Sigma(N_\Sigma)$.

Finally, for each connected component $\Sigma_i$ of $\Sigma$, let
$d_i$ be the number of elements in a generic stabilizer for the
$G$-action on the corresponding component $\tilde{P}_i$, and
$d_\Sigma:\,\Sigma\ra\N $ the function defined by the $d_i$'s.

For general values $\mu\in\Lambda_+$, let $\Sigma_\mu$,
$L_\mu$ etc. be defined by means of the shifting-trick.
The main result of this paper is the following

\bsat[Multiplicity Formula] \label{multf}
If $\mu\in\Lambda_+$ is a
regular value of $J$, the multiplicities $N^{(m)}(m\mu)$
are for $m>>0$ given by the formula
\be N^{(m)}(m\mu)=
            \int_{\Sigma_\mu} \frac{1}{d_{\Sigma_\mu}}\,
             \frac{Td\,({\Sigma_\mu})\,
             {Ch}^{\Sigma_\mu}(L^m_{\Sigma_\mu})}
             {D^{\Sigma_\mu}(N_{\Sigma_\mu})}.\label{orbifold1}\ee
If the image of the moment map, $J(M)$, is contained in $\G^*_{reg}$
(the set of regular elements of $\G^*$), one may set $m=1$ in this
formula:
\be N(\mu)= \int_{\Sigma_\mu} \frac{1}{d_{\Sigma_\mu}}\,\frac{{Td}\,
             ({\Sigma_\mu})\,
             {Ch}^{\Sigma_\mu}(L_{\Sigma_\mu})}
             {D^{\Sigma_\mu}(N_{\Sigma_\mu})},\label{orbifold}\ee
In particular, this is the case if $G$ is abelian.
\esat

The following sections are aimed at proving this Theorem.
We do not know whether the second part of the Theorem remains true
without the condition $J(M)\subset\G^*_{reg}$.\\

\noindent{\bf Remarks.}
\begin{enumerate}
\item
Comparing the right hand side of (\ref{orbifold}) to Kawasaki's
Riemann-Roch Formula for orbifolds \cite{K79}, the Theorem says that
$N(\mu)$ is equal to the Euler number of the orbifold-bundle
$L_\mu\ra M_\mu$. In particular, $N(\mu)$ is zero if the fiber of
$L_\mu$ over the smooth stratum of $M_\mu$ is a nontrivial
quotient of $\C$.

\item
Let $\Delta=J(M)\cap\T_+^*$, which is a convex polytope by
a result of Guillemin-Sternberg and Kirwan,
and $\Delta^*\subset
\Delta$ the set of regular values. By the Duistermaat-Heckman Theorem
\cite{DH82},
the diffeotype of the reduced space $M_\mu$ (and of course also of
$\Sigma_\mu$)
does not change as $\mu$ varies in a connected component of
${\rm int}(\Delta^*)$,
and the cohomology class of the symplectic form $\omega_\mu$
varies linearly.
In particular, the symplectic volume $\mbox{\rm Vol}(M_\mu)$ is
a polynomial on these connected components.
If the action of $G_\mu$ on $J^{-1}(\mu)$ is free, so that
${\Sigma_\mu}=M_\mu$, the right hand side  in (\ref{orbifold}) is equal to
a polynomial as well, since all that varies is the Chern character
${Ch}(L_\mu)=e^{\omega_\mu}$.
In the orbifold case, the behaviour is slightly more complicated:
For $\mu\in\Lambda_+$ in any given connected
component of ${\rm int}(\Delta^*)$, and any connected component
$\Sigma_{\mu,j}$ of $\Sigma_\mu$,
\be {Ch}^{\Sigma_\mu}(L_{\Sigma_\mu})|\,\Sigma_{\mu,j}
=\rho_\mu(g_j^{-1})\,c_L(g_j)\,e^{\omega_{\Sigma_\mu}},\ee
where $g_j\in G_\mu$ represents $\tau(\Sigma_{\mu,j})$, and
$\rho_\mu:G_\mu\ra S^1$
is defined by  $\rho_\mu(\exp(\xi))=e^{2\pi i\l \mu,\xi\r}$.
Hence, the right hand side of (\ref{orbifold}) is of the form
\be N(\mu)=\sum_j \rho_\mu(g_j^{-1})\, c_L(g_j)\,p_j(\mu),
\ee
where the $p_j$ are polynomials of degree
$\f{1}{2}\dim(\Sigma_{\mu,j})$.

\item
Since ${Ch}\,((L^m)_{m\mu})=e^{m\omega_\mu}$, the right hand side
of (\ref{orbifold1}) is a polynomial in $m$ if the $G_\mu$-action
on $J^{-1}(\mu)$ is free. In the orbifold-case, this is not true
in general since
\be {Ch}^{\Sigma_\mu}(L^m_{\Sigma_\mu})|\,\Sigma_{\mu,j}=
\rho_\mu(g_j^{-1})^m\,c_L(g_j)^m\,e^{m\omega_{\Sigma_\mu}}.\ee
\bdef[Ehrhart \cite{E77}]
A function $f:\N\ra \C$ is called an arithmetic polynomial,
if for some $k\in \N$, all the functions
\be q_j(m)=f(km-j),\,\,j=0,\ldots, k-1\ee
are polynomials. $k$ is called the period of $f$.
\ede
Equivalently, $f$ is an arithmetic polynomial
if and only if it can be written in the form
\be f(m)=\sum_{l=0}^{k-1}\exp(2\pi i\,\f{l\,m}{k})\,p_l(m),\ee
where the $p_l$ are polynomials. Taking $k$ such that $g^k=e$
for all $(g)\in \tau({\Sigma_\mu})$, the right hand side of (\ref{orbifold1})
clearly has this property.
\item
Our proof of Theorem \ref{multf} does not really use the
assumption that $M$ is K\"ahler. Everything will be derived
from the Equivariant Riemann-Roch Formula (\ref{char}),
which is of course valid in much more general situations.
Suppose for instance that $M$ is an arbitrary compact symplectic
manifold, equipped with a Hamiltonian $G$-action, and that these data
are quantizable. Then one can always choose a compatible, invariant
almost K\"ahler structure, and replace the virtual space
$\sum (-1)^i \,H^i(M,{\cal O}(L))$ with the index space of
some G-invariant Dirac operator for the Clifford module
$L\otimes \wedge(T^{(0,1)}M)^*$ (see \cite{BGV92,G94}).
As an immediate consequence of the Berline-Vergne Formula for
the character, Theorem \ref{RiemannRoch} below,
the multiplicities $N(\mu)$ defined in this way
do not depend on the choice of the almost K\"ahler structure or
of the quantizing line bundle $L$.
\end{enumerate}

\bigskip

\noindent{\bf Example:}
Let $M=\C P(2)$, equipped with the Fubini-Study K\"ahler form
$\omega_{FS}$.
Let $G=S^1$ act by
\[e^{i\phi}.[z_0:z_1:z_2]=[e^{i\phi}z_0:e^{-i\phi}z_1:z_2].\]
This action is Hamiltonian, and has a moment map
\[J([z_0:z_1:z_2])=
\frac{|z_1|^2-|z_0|^2}{|z_0|^2+|z_1|^2+|z_2|^2}.\]
The dual of the tautological line bundle serves as a quantizing
line bundle $L$.
We also consider the tensor powers $L^m$, which are
quantizing line bundles for $(M,m\,\omega_{FS})$.
By Kodaira's Theorem, $H^i(M,\O(L^m))=0$ for all $i>0,\,m\in\N$.
If we identify the
spaces $\Gamma_{hol}(M,L^m)$ with the homogeneous polynomials of
degree $m$ on $\C^3$, the representation of $S^1$ is induced
by the action $e^{i\phi}.(z_0,z_1,z_2)=
(e^{i\phi}z_0,e^{-i\phi}z_1,z_2)$ on ${\C^3}$.
The isotypical subspace of $\Gamma_{hol}(M,L^m)$
corresponding to the weight $l\in\Z$
is, for $l\ge 0$, spanned by
\[ z_0^l\,z_2^{m-l},\,z_0^{l+1}\,z_1\,z_2^{m-l-2},\ldots,
z_0^{l+r}\,z_1^r\,z_2^{m-l-2r}\]
with $r=\left[\frac{m-l}{2}\right]$. For $l \le 0$, the roles
of $z_0$ and $z_1$ are reversed. Thus
\[ N^{(m)}(l)=\left\{\begin{array}{ll}
              {1+\left[\f{m-|l|}{2}\right]}&{\mbox{ if $|l|\le m$}}\\
                       {0}&{\mbox{ otherwise}}
               \end{array}  \right.,\]
for all $l\in\Z,\,m\in\N$.
On the other hand,
the image of the moment map $J^{(m)}=mJ$ is
the interval $-m\,\le\, \mu\,\le\, m$, with critical values at
$-m,\,0,\,m$. If $0< \, |l|\,< m $, the level set $(mJ)^{-1}(l)$
consists of two orbit type strata: On the set where $z_2\not=0$,
the action is free, and on the set where $z_2=0$, the stabilizer
is $\Z_2$. Writing $S=\{e,g\}$, the reduced space $M^e_\mu=M_\mu$
is an orbifold with a $\Z_2$ singularity (the
``teardrop-orbifold''), and $M^g_\mu$ is the singular point.

Since $c_{L^m}(g)=(-1)^m$ and $\rho_l(g^{-1})=(-1)^l$, we expect
the multiplicities to grow like
$N^{(m)}(l)=p_e(m,l)+(-1)^{(m-l)}p_g(m,l)$ where $p_e$ is a first
order polynomial and $p_g$ a constant. Comparison with the
above formula shows that this is indeed the case, with
$p_e(m,l)=3/4+(m-|l|)/2$ and $p_g(m,l)=1/4$.
Note that for $m-|l|$ even, the fiber of $L^m_l$ at the singular
point is $\C$, whereas form $m-|l|$ odd it is $\C/\Z_2$. This means
that for $m-|l|$ odd, all holomorphic sections of $L^m_l$ have to
vanish at the singular point. Again, this fits with the above
explicit formulas.

\section{Some equivariant cohomology}\label{section3}
We start by reviewing Cartan's model for equivariant cohomology,
following Berline and Vergne \cite{BV85}.
Let $M$ be a compact manifold, $G$ a compact
Lie group, and $\Phi:G\times M\ra M$
a smooth action. Denote by ${\cal A}_G(M)$ the space of $G$-invariant
polynomial mappings $\sigma:\G\ra \A(M)$, that is, $\sigma(\xi)$ depends
polynomially on $\xi$ and satisfies the equivariance property

\be \sigma(\mbox{Ad}_g(\xi))=\Phi_{g^{-1}}^*(\sigma(\xi)).           \ee

The elements of $\A_G(M)$ are called {\em equivariant differential forms},
and the space ${\cal A}_G(M)$ is preserved under the
{\em equivariant differential}
\be \D:\A_G(M)\ra\A_G(M),\,\,(\D\sigma)(\xi)=
\d(\sigma(\xi))+2\pi i \big(\iota(\xi_M)\sigma (\xi)\big).\ee
Here,
$\xi_M$ denotes the fundamental vector field, i.e. the generating vector
field of the flow $(t,p)\mapsto \exp(-t\xi).p$.

Equivariance together with Cartan's identity for the Lie derivative,
$\L_Y=\iota_Y\circ \d +\d\circ \iota_Y$, implies $\d_\G^2=0$.
The cohomology $H_G(M)$ of the complex
$({\cal A}_G(M),\D)$ is called the {\em equivariant cohomology}.
One can show that if the action of $G$ on $M$
is locally free, the pullback mapping ${\cal A}(M/G)\ra {\cal A}(M)^G_{hor}
\hra {\cal A}_G(M)$ gives rise to
an isomorphism  $H(M/G)\ra H_G(M)$.
After choosing a principal connection
on $M\ra M/G$, the inverse is induced on the level of forms by the
mapping
\be {\cal A}_G(M)\ra {\cal A}(M)^G\ra {\cal A}(M)^G_{hor}\cong {\cal A}
(M/G)\ee
given by substituting $\f{i}{2\pi}$ times
the curvature in the $\G$-slot, followed by
projection onto the horizontal part (for the proof, see \cite{DV93}).

For what follows, it will be necessary to relax the polynomial
dependence on $\xi$ to analytic dependence, possibly defined only
on some neighborhood of $0\in\G$.
We will denote the corresponding space of equivariant forms by
${\cal A}^\omega_G(M)$, and the
cohomology by $H_G^\omega(M)$.

Suppose now that ${\cal V}\ra M$ is a $G$-equivariant Hermitian
vector bundle over $M$, with fiber dimension $N$. Let ${\cal A}(M,{\cal V})$
be the bundle-valued differential forms, and ${\cal A}_G(M,{\cal V})$
their equivariant counterpart.
For each $G$-invariant Hermitian connection $\nabla:{\cal A}(M,{\cal V})\ra
{\cal A}(M,{\cal V})$,
the moment map $\mu\in{\cal A}_G(M,\mbox{End}({\cal V}))$ of Berline
and Vergne is defined by
\be \mu(\xi).\sigma :=
\xi.\sigma - \nabla_{\xi_M}\sigma,\ee
where $\sigma\ra \xi.\sigma$ denotes the representation
of $\G$ on the space of sections.
Geometrically, $\mu(\xi)$ is the vertical part (with respect to the
connection) of the fundamental
vector field $\xi_\V$ on $\V$. Let $F({\cal V})\in{\cal A}^2
(M,\mbox{End}({\cal V}))$ denote the curvature of $\nabla$.
The {\em equivariant} curvature is then defined by
\be F_\G({\cal V},\xi)=F({\cal V})+2\pi i \mu(\xi),\ee
and it satisfies the Bianchi identity with respect to the
equivariant covariant derivative $\nabla_\G=\nabla+2\pi i \iota(\xi_M)$.
Suppose now that $A\ra f(A)$ is the germ of a $U(N)$-invariant
analytic function on ${\frak{u}}(N)$. Then $f(F_\G)\in{\cal A}_G(M)$
is $\d_\G$-closed, and one can show that choosing a different
connection changes $f(F_\G)$ by a $\d_\G$-exact form.
The corresponding cohomology classes are called the {\em equivariant
characteristic classes} of ${\cal V}\ra M$.
If the action on $M$ is
locally free, one can choose $\nabla$ in such a way that $\mu=0$,
which shows that
the mapping $H^\omega_G(M)\ra H(M/G)$ sends the
equivariant characteristic classes of ${\cal V}$ to the usual
characteristic classes of the orbifold-bundle ${\cal V}/G$.

The following
characteristic classes will play a role in the sequel:

\begin{itemize}
\item[(a)] The equivariant Chern character, defined by
            \be {Ch}_\G(\V,\xi)=\mbox{tr}
(e^{\f{i}{2\pi}F_\G(\V,\xi)}).\ee
In the above geometric quantization setting, ${\cal V}=L$ is a line
bundle, and for the equivariant curvature one has
$\f{i}{2\pi}F_\G(\V,\xi)=\omega+2\pi i\l J,\xi\r $, thus
\be Ch_\G(L,\xi)=e^{\omega+2\pi i\l J,\xi\r}.\ee
More generally, if  $g\in G$ acts trivially on the base $M$,
one defines
\be {Ch}_\G^g(\V,\xi)=\mbox{tr}(g^\V\,\,e^{
\f{i}{2\pi}F_\G(\V,\xi)})\ee
where $g^\V\in\Gamma(M,\mbox{End}({\cal V}))$
is the induced action of $g$.
In the line bundle case,
this is simply $c_L(g){Ch}_\G(L,\xi)$ where
where $c_L(g)\in S^1$ is the action of $g$ on the fibers.
\item[(b)]  The equivariant Todd class,
           \be {Td}_\G({\cal V},\xi)=
           \det\Big(\f{\f{i}{2\pi}F_\G(\V,\xi)}{1-e^{-\f{i}{2\pi}
           F_\G(\V,\xi)}}\Big).\ee
           The Todd class of a complex manifold is defined as the Todd class
           of its tangent bundle.
\item[(c)] The equivariant Euler class
           \be\chi_\G({\cal V},\xi)=\det(\f{i}{2\pi}F_\G(\V,\xi)).\ee
\item[(d)] The class
          \be D_\G^g(\V,\xi)=\det (I-(g^{-1})^\V e^{-\f{i}{2\pi}
           F_\G(\V,\xi)}),\ee
           for $g\in G$ acting trivially on $M$.
\end{itemize}

All of this also makes sense for symplectic vector bundles, since the
choice of a compatible complex structure reduces the structure group
to $U(n)$, and any two such choices are homotopic.

The equivariant Euler class occurs in the Localization Formula of
Atiyah-Bott and Berline-Vergne.

\bsat[Atiyah-Bott \cite{AB84}, Berline-Vergne \cite{BV83}] \label{AB84}
Suppose $M$ is an orientable $G$-ma\-ni\-fold, and
$\sigma\in {\cal A}_G^\omega(M)$ is
$\D$-closed. Assume that $\xi\in\G$
is in the domain of definition of $\sigma(\xi)$, and
let $F$ be the set of zeroes of $\xi_M$. The connected components of
$F$ are then
smooth submanifolds of even codimension, and the normal
bundle $N_F$ admits a Hermitian structure which is invariant under
the flow of $\xi_M$.
Choose orientations on $F$ and $M$ which are compatible with the
corresponding orientation of $N_F$. Then
\be \int_M \sigma(\xi)=\int_{F}
\frac{\iota_F^*\sigma(\xi)}{\chi_\G(N_F,\xi)},\label{localization} \ee
where $\iota_F:F\ra M$ denotes the embedding.
\esat

We will need this result only in the symplectic
or complex case, where the above orientations are given in a natural way.

Jeffrey and Kirwan \cite{JK93} have proved a different sort
of Localization Formula for the case of {\em Hamiltonian} G-spaces.
We will need a stationary phase version of their result.
Let $(M,\omega)$ be a Hamiltonian $G$-space, with moment map $J:M\ra\G^*$,
and suppose $0$ is a regular value of $J$. Let
$\sigma\in{\cal A}^\omega_G(M)$ be $\D$-closed,
and let $\Delta\in C^\infty_0(\G)$ be a cutoff-function, with $\sigma(\xi)$
defined for $\xi\in\mbox{supp}(\Delta)$, and $\Delta=1$ in a neighborhood of
$0$.
Consider the integral
\be\int_\G\int_M  \Delta(\xi)\sigma(\xi)e^{m(\omega+2 \pi i\,\l J,\xi\r)}
\,\d\xi,\ee
where $\d\xi$ is the measure on $\G$ corresponding to the normalized measure
on $G$. Notice that the ``equivariant symplectic form''
$\omega+2 \pi i\,\l J,\xi\r$ is $\D$-closed.
Since $e^{m\omega}$ is simply a polynomial in $m$, the
leading behaviour of this integral for $m\ra\infty$ is determined by the
stationary points of the phase function $e^{2\pi i m\l J,\xi\r}$.
Stationarity in $\G$-direction gives the condition $J=0$, and stationarity
in $M$-direction
the condition $\d \l J,\xi\r=0$, or $\xi=0$ since the action on $J^{-1}(0)$
is locally free.
Let $M_0= J^{-1}(0)/G$ be the reduced space, $\pi:J^{-1}(0)\ra
M_0$ the projection and $\iota:J^{-1}(0)\ra M$ the embedding.
Consider the mapping
\be \kappa:H^\omega_G(M)\ra H(M_0),\label{kappa}\ee
given by composing pullback to $J^{-1}(0)$ with the mapping
$H^\omega_G(J^{-1}(0))\ra H(M_0)$. On the level of forms
$\sigma\in {\cal A}^\omega_G(M)$, the form $\pi^*\kappa(\sigma)$
is by definition equal to the horizontal part of
$\iota^*\sigma(\f{i}{2\pi}F^\theta)$, where
$F^\theta\in{\cal A}^2(J^{-1}(0),\G)^G_{hor}$
is the curvature of some connection
$\theta\in{\cal A}^1(J^{-1}(0),\G)^G$.

\bsat \label{jf}
For $m\ra\infty$,
\be \int_\G\int_M
\Delta(\xi)\sigma(\xi)e^{m(\omega+2\pi i\l J,\xi\r)}\,\d\xi
=\f{1}{d}\int_{M_0}\kappa(\sigma) \,
e^{m\omega_0}\,\,\,+O(m^{-\infty}),\label{jklocf} \ee
where $d$ is the number of elements in the generic stabilizer
for the $G$-action on $J^{-1}(0)$.
\esat

\bbew

This is a variation of Theorem 4.1 in  Jeffrey-Kirwan
\cite{JK93}, which deals with polynomial equivariant
cohomology classes, and where one has a Gaussian convergence
factor instead of the cutoff.

Following \cite{JK93}, we will perform the integral in a local
model for $M$ near $J^{-1}(0)$.
By the Coisotropic Embedding Theorem of Gotay, a neighborhood of
$J^{-1}(0)$ in $M$ is equivariantly symplectomorphic to
a neighborhood of the zero section of the trivial bundle
$J^{-1}(0)\times \G^*$,
with symplectic form $\pi^*\omega_0+\d\l \alpha,\theta\r$, where
$\alpha$ is the coordinate function on $\G^*$. In this model,
$G$ acts by $g.(x,\alpha)=(g.x,\mbox{Ad}_{g^{-1}}^*(\alpha))$ ,
and the moment map is
simply $J(x,\alpha)=\alpha$. Using the model and another
cutoff-function $\Delta'(\alpha)$ on $\G^*$, equal to $1$ near the origin
and with sufficiently small support, the same computation
as in \cite{JK93} shows that the integral is equal to
\[ m^{\dim G}\int_{J^{-1}(0)}
\int_{\G^*}\int_{\G} \Delta(\xi)\Delta'(\alpha)
(\iota^*\sigma)(\xi)e^{m(\pi^*\omega_0+\l \alpha,F^\theta+2\pi i
\xi\r)}\, \d\xi\, \d\alpha\,\d g +O(m^{-\infty}).\]
Here, $\d g$ denotes the (vertical) volume form on the fibers of
$J^{-1}(0)\ra M_0$, corresponding to the canonical identification
$T_x(\mbox{fiber})\cong\G$ by means of the $G$-action.

Now apply the Stationary Phase Theorem (see e.g. \cite{H90}, Theorem 7.7.5)
to the $\alpha,\xi$-integral, the relevant phase function
being $e^{2\pi i m \l \alpha,\xi\r}$.
Since $e^{m\l \alpha,F^\theta\r}$
is simply a polynomial in $\alpha$, the stationary phase
expansion
terminates after finitely many terms, and the result is
\[ \int_{J^{-1}(0)} e^{m \pi^*\omega_0}
\sum_{r=0}^\infty \f{1}{r!}\big(\f{i}{2 \pi m}\big)^r\left.
\Big(\sum_j\f{\p}{\p \xi_j}
\f{\p}{\p \alpha^j}\Big)^r\right|_{\atop\stackrel{\xi=0}{\alpha=0}}
\iota^*\sigma(\xi)\, e^{m\l \alpha,F^\theta\r }
\,\,\d g +O(m^{-\infty})\]
\[ =\int_{J^{-1}(0)} \pi^*\,(e^{m \omega_0})\,\iota^*\sigma(\f{i}{2\pi}
F^\theta)\,\d g +O(m^{-\infty}).\]
Since $\iota^*\sigma(\f{i}{2\pi}F^\theta)$ gets wedged with $\d g$,
only its horizontal part, which by definition of $\kappa$
is  $\pi^*\kappa(\sigma)$,
contributes to
the integral. The result (\ref{jklocf}) now follows by
integration over the fiber; the factor
${1}/{d}$ appears since this is the volume of a generic fiber.
\ebew

We now turn to the Equivariant Hirzebruch-Riemann-Roch Theorem, in
the form due to Berline and Vergne \cite{BV85}. Let $M$ be a compact
complex manifold, equipped with a holomorphic action of $G$,
and let $L\ra M$ be a G-equivariant holomorphic line bundle.
Define the character $\chi\in R(G)$ as in (\ref{char}).
\bsat\label{RiemannRoch}
For  $\xi$ sufficiently close to zero,
\be \chi(e^{\xi})=
\int_M {Td}_\G(M,\xi)\,{Ch}_\G(L,\xi).\label{e}\ee
More generally, if  $g\in G$, one has for all sufficiently small
$\xi\in \K$, the Lie algebra of the centralizer $Z_g$ of $g$:
\be \chi(g\, e^{\xi})=
\int_{M^g}
\f{{Td}_\K(M^g,\xi)\,{Ch}_\K^g
(L|\,M^g,\xi)}{ D^g_\K(N^g,\xi)},
\label{g}\ee
where $M^g$ is the  fixed point set and  $N^g\ra M^g$ its normal bundle.
\esat

To be precise, Berline and Vergne have shown how to rewrite the Equivariant
Atiyah-Segal-Singer Index Theorem for Dirac operators in this style,
with the
equivariant $\hat{A}$-genus appearing
on the right hand side. This formula, however,
implies the above Theorem
in the same way as the usual Atiyah-Singer Index Theorem
leads to the
Hirzebruch-Riemann-Roch Formula; see \cite{BGV92}, p. 152 for the
calculations.

\section{The stationary phase approximation}

In this section, we will prove the first part of Theorem \ref{multf}.

By the shifting trick, it is sufficient to consider
the case $\mu=0$.
The idea is to substitute the expressions from the Equivariant
Hirzebruch-Riemann-Roch Theorem \ref{RiemannRoch} for $\chi^{(m)}$ in
\be N^{(m)}(0)=\int_G \chi^{(m)}(h)\d h,\label{4.5}\ee
and apply the Localization Formula, Theorem \ref{jf}.
For this, we need to know what happens to the equivariant
Todd class of $M$ under the mapping (\ref{kappa}):

\blem\label{todd} Let
\be j_\G(\xi)=\det\Big({\ts\f{1-e^{-{ad}(\xi)} }
    {{ad}(\xi)}}\Big)\ee
be the Jacobian of the exponential mapping $\exp:\G\ra G$.
Then
\be \kappa({Td}_\G(M) j_\G)={Td}\,(M_0).\ee
\elem

\bbew
Identify the vertical subbundle of $TJ^{-1}(0)$ with the trivial
bundle $\G$, and the symplectic bundle $\G\oplus I\G$ (where $I$
is the complex structure of $M$) with $\G_\C$. Then
\[ \iota^*(TM)=\pi^*(TM_0)\oplus \G_\C.      \]
Since the equivariant Todd class of $\G_\C$ is just $j_\G^{-1}$,
this shows $\kappa({Td}_\G(M))=Td\,(M_0)\,\kappa(j_\G^{-1})$, q.e.d.
\ebew

We will now consider the contribution to (\ref{4.5}) coming from
a small ${\rm Ad}$-invariant
neighborhood of a given orbit $(g)={\rm Ad}(G).g$.
Let $\sigma\in C^\infty(G)$ be an
$\mbox{Ad}$-invariant cutoff-function, supported in a sufficiently
small neighborhood of $(g)$ and equal to 1 near $(g)$.
Consider the integral
\be I_g(m)=\int_G \sigma(h)\chi^{(m)}(h)\d h.\ee
Since (\ref{g}) only holds for
$\xi\in\K$, we want to replace this integral by an integral over $Z_g$.
(Of course, this step is void in the abelian case.)

Let ${\frak r}\subset\G$
be the orthogonal of $\K$ with respect
to some invariant inner product (or, more intrinsically, the
annihilator of $(\G^*)^g\cong\K^*$). For $k\in Z_g$, let
$k^{\frak{r}}$ denote the action of $k$ on ${\frak r}$.

\blem \label{lemmaGK}
Let $f\in C^\infty(G)$ be $\mbox{Ad}(G)$-invariant, with support
in a small neighborhood of $\mbox{Ad}(G).g$. Then, for a suitable
$\mbox{Ad}(Z_g)$-invariant cutoff-function
$\tilde{\sigma}\in C^\infty_0(Z_g)$,
supported near $e\in Z_g$ and identically $1$ near $e$,
\be \int_G f(h)\d h=\int_{Z_g} f(gk)
    \det(I-(g^{-1}k^{-1})^{\frak{r}})\tilde{\sigma}(k)\,\d k.
\label{GK}\ee
\elem

\bbew
 From the Slice Theorem for actions of compact Lie groups, it
follows that an invariant neighborhood of the orbit
$\mbox{Ad}(G).g$ is equivariantly diffeomorphic to a neighborhood
of the zero section of $G\times_{Z_g}\K$, where $Z_g$ acts on $\K$
by the adjoint action. Using the exponential map for $Z_g$,
it follows that the mapping
\be \phi:\,G\times_{Z_g}Z_g\ra G,\,\,(h,k)\mapsto \mbox{Ad}_h(gk), \ee
with $Z_g$ acting on itself by $\mbox{Ad}$,
is an equivariant diffeomorphism from a neighborhood of the unit
section to a neighborhood of the orbit.
Let $\d\nu$ be
the canonical measure on the group bundle $G\times_{Z_g}Z_g=:W$,
constructed
from the normalized invariant measures on $Z_g$ and $G/Z_g$.
We have to compute the tangent mapping to $\phi$, but
by equivariance it is sufficient to do this at points
$(e,k)$ in the fiber over $e\,Z_g$, which is canonically
isomorphic to $Z_g$. If we identify $T_{(e,k)}W\cong
\frak{r}\oplus\frak{k}$, and $TG\cong G\times\G$ using
left trivialization, the tangent mapping is given by
\[ T_{(e,k)}\phi(\xi,\eta)=((1-\mbox{Ad}_{(gk)^{-1}})(\xi),
\eta).\]
This shows that the measure transforms according to
\be \phi^*\,\d g = \det(I-(g^{-1}k^{-1})^{\frak r})\,\d \nu.\ee
We can now perform the integral by first pulling $f$ back to
$G\times_{Z_g}Z_g$, multiplying with a suitable cuttoff function
$\tilde{\sigma}$,
integrating over the fibers of $W\ra G/Z_g$, and then
integrate over the base $\mbox{Ad}(G).g\cong G/Z_g$.
\ebew

Using the Lemma, we find that
\be I_g(m)=
\int_{Z_g} \tilde{\sigma}(k)\chi^{(m)}(gk)
\det(I-(g^{-1}\,k^{-1})^{\frak r})\,
\d k.\ee
Replacing this with an integral over the Lie algebra $\K$, and using
(\ref{g}) gives
\be I_g(m)=\int_\K\int_{M^g} \Delta(\xi)
\f{{Td}_\K(M^g,\xi){Ch}_\K^g(L^m|M^g,\xi)}{ D^g_\K(N^g,\xi)}
\det(I-(g^{-1}e^{-\xi})^{\frak r})
j_\K(\xi)\d \xi,\label{LK}\ee
with $\Delta(\xi)=\tilde{\sigma}(e^{\xi})$.
Let $\kappa_g:\,H_{Z_g}(M^g)\ra H(\Sigma_g)$ be the mapping defined
by (\ref{kappa}), with $M$ replaced by $M^g$ and $G$ by $Z_g$.
By Lemma \ref{todd},
\be \kappa_g(Ch^g_\K(L^m|M^g))\,\kappa_g(j_\K)=
Ch^{\Sigma}(L^m_\Sigma)|\Sigma_g.\ee

For $x\in P^g=M^g\cap J^{-1}(0)$, let
${\frak r}_M(x):=\{\xi_M(x)|\xi\in {\frak r}\}\cong {\frak r}$.
Then
\be N^g(x)=N_\Sigma(y)\oplus{\frak r}_M(x)\oplus I{\frak r}_M(x)
=N_\Sigma(y)\oplus {\frak r}_\C,\ee
where $y=G.(x,g)$.
But $D^g_\K({\frak r}_\C,\xi)=\det(I-(g^{-1}e^{-\xi}))$, hence
\be \kappa_g(D_\K^g(N^g,\xi))=
\kappa_g(\det(I-(g^{-1}e^{-\xi})^{\frak r}))\,\,
D^\Sigma(N_\Sigma)|\Sigma_g.\ee
With these preperations, we apply Theorem \ref{jf} to the
integral \ref{LK}, and obtain
\be I_g(m)={\sum}'_{j} \frac{1}{d_j}
            \int_{\Sigma_j} \frac{Td\,(\Sigma_j)
             {Ch}^\Sigma(L^m_\Sigma)}
             {D^\Sigma(N_\Sigma)}+O(m^{-\infty}),\ee
the sum being over the connected components of $\Sigma_g$.
Summing over all contributions, we get (\ref{orbifold1}) up to an error
term $O(m^{-\infty})$.
As we remarked above, the right hand side of (\ref{orbifold1})
is an arithmetic polynomial as a function of $m$.
But if $f:\N\ra \Z$ is an  integer-valued function
with $\lim_{m\ra\infty}(f(m)-p(m))=0$ for some polynomial $p$,
then $f(m)=p(m)$ for large $m$.
This shows that the error term is zero for large $m$, and finishes
the proof of the first part of Theorem \ref{multf}.

\section{Counting lattice points}

To prove the second part of Theorem \ref{multf}, i.e. that we can
set $m=1$, all we have to show is that $m\ra N^{(m)}(m\mu)$ is
an arithmetic polynomial.

\bsat \label{stepw}
Suppose that $J(M)\subset \G^*_{reg}$. Then the function
$m\mapsto N^{(m)}(m\mu)$ is an arithmetic polynomial for all
$\mu\in\Lambda_+$.
\esat

Before we prove this, we convert the computation of the multiplicities
into a problem of counting lattice points.
The next steps are based on work of
Guillemin-Lerman-Sternberg
\cite{GLS88} and Guillemin-Prato \cite{GP90},
except that we replace their use of the
Atiyah-Bott Lefschetz Formula with the Localization Formula \ref{AB84}
applied to (\ref{e}), since we do not want to assume
isolated fixed points.
Consider the action of the maximal torus $T\subset G$, with
its moment map $J^T$ equal to $J$ followed by projection
to $\T^*$.

\bprop
\label{prop5.1}
For all generic $\xi\in\T$,
\be \chi(e^{\xi})=\sum_{\F} \int_\F \frac{{Td} (\F)e^{\omega
+ 2\pi i \l J_\F,\xi\r}}
{\det(I-e^{-\f{i}{2\pi}F_\T(N_\F,\xi)})},\label{fixedpoints}\ee
the sum being over the fixed points manifolds of the
$T\subset G$-action, $N_\F$ the corresponding normal bundles,
and $J_\F$ the constant value of $J$ on $\F$.
\eprop

\bbew
By ``generic'' we mean that the zero set of
$\xi_M$ is equal to the fixed point set of the $T$-action.
To get ({\ref{fixedpoints}), apply the
Localization Formula, Theorem \ref{AB84}, to (\ref{e}).
The bundle $TM|_\F$ splits
into the direct sum $T\F\oplus N_\F$. Since $T$ does not act on
$T\F$,
\[{Td}_\T(M,\xi)={Td}\,(\F)\,{Td}_\T(N_\F,\xi).\]
The equivariant Euler class of $N_\F$ in Proposition
\ref{AB84} cancels
the denominator of the equivariant Todd class, which immediately
gives (\ref{fixedpoints}).
\ebew

Although the left hand side of (\ref{fixedpoints})
is an analytic function of $\xi$, the
individual summands on the right hand side have poles.
Since they are not in $L^1_{loc}$, they do not a priori
define distributions
on $\T$.
This problem can be fixed as follows \cite{D93}.
By using the splitting principle (or simply a partition of unity
on ${\cal F}$) if necessary, we can assume
that  $N_\F$ splits into a direct sum of invariant line bundles
$N_\F^1,\ldots,N_\F^r$. Let
$\alpha^j_\F$ be the weight for the
$T$-action on $N_\F^j$, that is, $e^\xi\in T$ acts by the
character $\exp(2\pi i \l \alpha_\F^j,\xi\r)$.
Each $\alpha^j_\F$ determines an orthogonal hyperplane in $\T$,
let $C$ be any fixed connected component in the complement of
the union of all these hyperplanes.
If we replace $\xi$ by $\xi-i\eta$ in
({\ref{fixedpoints}), with $\eta\in C$, the terms on the right hand side
are analytic for all $\xi$. One can therefore regard
(\ref{fixedpoints}) as an equality of distributions, with the
summands on the right hand side defined defined as a
distributional limit for $\eta\ra 0$ in $C$.

Let us now first discuss the abelian case, i.e. assume that
$G=T$ is a torus. Denote by
$F^j(N_\F)$ the components of the curvature.
By expanding
$\det(I-e^{-\f{i}{2\pi}{F}_\T(N_\F,\xi-i\eta)})^{-1}$
into its Taylor series w.r.t. ${F}^j(N_\F)$, we can write it
as a finite sum
\be \det(I-e^{-\f{i}{2\pi}{F}_\T(N,\xi-i\eta)})^{-1}=
    \sum_{s\in\N^r} \f{p_s({F}^1(N_\F),\ldots,
    {F}^r(N_\F))}
    {\prod_j (1-e^{-2\pi i\l {\alpha}^j_\F,\xi-i\eta\r})^{s_j}},
\label{5}
\ee
where for all $s=(s_1,\ldots,s_r)$, $p_s$ is a polynomial.
We now invoke the ``polarization trick'' used  in ref.
\cite{GLS88,GP90}.
For each $\alpha^j_\F$,
write
\be \check{\alpha}^j_\F=\left\{\begin{array}{ll}
{\alpha}^j_\F&\mbox{ if }
\l {\alpha}^j_\F,\eta\r>0\\
-{\alpha}^j_\F&\mbox{ if }
\l {\alpha}^j_\F,\eta\r<0\end{array}\right.
\label{signs}\ee
for any, hence all, $\eta\in C$.
Let $l_j^0=0$ if $\check{\alpha}^j_\F=\alpha^j_\F$, $1$ otherwise.
Then
\be \chi(e^{\xi-i\eta})=\sum_\F\sum_{s\in \N^r}c_{\F,s}
\f{e^{2\pi i\l J_\F-\sum l_j^0 s_j\check{\alpha}^j_\F,\xi-i\eta\r}}
{\prod_j (1-e^{-2\pi i\l \check{\alpha}^j_\F,\xi-i\eta\r})^{s_j}}
\label{lhs}\ee
with
\be c_{\F,s}= (-1)^{k_{\F,s}} \int_\F {Td}\,(\F)\,e^\omega p_s({F}^1
(N_\F),\ldots,
{F}^r(N_\F)),\ee
where $k_{{\cal F},s}=\sum l_j^0 s_j$ is the number of sign changes.

For given $\F,s$, write $(a^1,\ldots,a^N)$
for the list of $\check{\alpha}^j_\F$'s, appearing with respective
multiplicities $s_j$. Since
\be (1-e^{-2\pi i\l a^j,\xi-i\eta\r})^{-1}=\sum_{l_j=0}^\infty
e^{- 2\pi i \l l_j\,a^j,\xi-i\eta \r},\ee
we get
\be \chi(e^\xi)=\sum_{\F,s} c_{\F,s}\sum_{l\in\Z^N_+}
   e^{2\pi i\l J_\F-\sum (l_j+l_j^0)a^j,\xi\r}\ee
(the sum over $\Z^N_+:=\{l\in \Z^N:\,l_j\ge 0\}$
is a well-defined periodic distribution).
Comparing this to
\be \chi(e^\xi)=\sum_{\mu\in\Lambda}
N(\mu)e^{2\pi i\l\mu,\xi\r}\label{abelian}\ee
yields
\be N(\mu) = \sum_{\F,s}c_{{\cal F},s}{\frak P}_{\F,s}(J_\F-\mu-
\sum l_j^0 a^j)\label{latticepoints}\ee
where the partition function ${\frak P}_{\F,s}(\nu)$
is the number of solutions
$k\in \Z^N$ of $\sum k_j a^j=\nu$, $k_j\ge 0$.

Starting from this expression, we will now show that $N^{(m)}(m\mu)$
is an arithmetic polynomial. We have to replace $\omega$ by
$m\omega$, $\mu$ by $m\mu$ and $J$ by $mJ$.
Since
$c_{\F,s}^{(m)}$ is a polynomial in $m$, it is sufficient to show that
the number of integer solutions of
\be m(J_\F-\mu)=\sum (l_j+l_j^0) a^j,\,\,l_j\ge 0 \label{10}\ee
is an arithmetic polynomial as a function of $m$.
Let us write $\nu=J_\F-\mu$, and consider $A=(a^1,\ldots,a^N)$ as a
$\Z$-linear mapping $\Z^N\ra \Z^p$, where $p=\dim(T)$. We are thus
looking for integer solutions of
\be m\nu= A\,l,\,\,l_j\ge l_j^0.\label{11}\ee
We will need the following

\bsat[Ehrhart \cite{E77}]
Let $L$ be a lattice,  with underlying vector space $L_\R=L\otimes_\Z \R$,
and $\Delta\subset L_\R$ a lattice polytope,
i.e. a polytope whose vertices are all lattice points. Then,
for all $r\in \N$, the counting function
 \be f(m)=\#\big(\f{m}{r}\Delta\cap L\big) \ee
is an arithmetic polynomial, with period $r$.
\esat

Let now $x_0\in\R^N$ be any solution of $Ax=\nu$. The general solution
of $A\,x=m\nu$ is thus given by the affine plane
$E_m=mx_0+\mbox{ker}(A)$. Let $r\in \N$ be the smallest number such
that the vertices of the polytope
$\Delta:=E_r\cap \R^N_+$ are lattice points. If $l^0=0$, the set of
solutions of (\ref{11}) is the intersection $\f{m}{r}\Delta\cap \Z^N$,
so the number of solutions is an arithmetic polynomial by Ehrhart's Theorem.
If $l^0\not=0$, let $\Delta_j$ be the face of $\Delta$ defined by
$x_j=0$, and let $\Delta'$ be the union of all $\Delta_j$ for which
$l_j^0=1$. Then the solution set of (\ref{11}) is
\[ \f{m}{r}\Delta\cap \Z^N\,-\,\f{m}{r}\Delta'\cap \Z^N,\]
and this is again an arithmetic polynomial by Ehrhart's Theorem.
This proves Theorem \ref{stepw} in the abelian case.

Suppose now that $G$ is nonabelian, but that $J(M)$ is contained
in the set of regular elements, $\G^*_{reg}=G.\mbox{int}(\T^*_+)$.
We will show how this reduces to the abelian case.
By the Symplectic Slice Theorem \cite{GS84}, $Y_+=J^{-1}(\mbox{int}
(\T^*_+))$ is a symplectic (but not necessarily K\"ahler)
submanifold of $M$, and is in fact a
Hamiltonian $T$-space, with the restriction of $J$ serving
as a moment map.
The above assumption implies that $Y_+$ is a {\em closed} submanifold, and
$M=G\times_T Y_+$. The restriction $L_+=L|Y_+$ renders a quantizing
bundle for $Y_+$. Consider the expression
\be \chi'(e^\xi):= \int_{Y^+} {Td}_\T(Y_+,\xi)\,
    {Ch}_\T(L_+,\xi).\ee
We claim that this is of the form
\be \chi'(e^\xi)=\sum_{\mu\in\Lambda} N'(\mu) \,e^{2\pi i \l \mu,\xi\r},\ee
where $N'(\mu)\not=0$ for only finitely many lattice points, and
$N'(\mu)=0$ unless $\mu\in J(Y_+)\cap\Lambda\subset\Lambda_+$.
Indeed, one can check directly that $\chi(e^\xi)$ comes from
a function on $T$, given near any point $g\in T$
by the formula (\ref{g}), and then repeat
the above analysis. (One can also pick a
$T$-invariant almost K\"ahler structure on $M$, and then realize
$\chi(e^\xi)$ as the equivariant index for some Dirac operator
associated to the Clifford module $L\otimes
\Lambda(T^{(0,1)}Y_+)^*$.)

\blem For all $\mu\in\Lambda_+$, $N(\mu)=N'(\mu)$.
\elem

Since we know that ${N'}^{(m)}(m\mu)$ is an arithmetic polynomial, this
will finish the proof of Theorem \ref{multf}.\\

\bbew
Let us go back to the formula (\ref{fixedpoints}) for the character.
Notice that the Weyl group $W=N_G(T)/T$ acts on $M^T$ by
permuting the connected components, and that $M^T$ consists of
its portion in $Y_+$ and the $W$-transforms thereof.
Let ${\cal F}\subset Y_+$ be a connected component of $M^T$. The
normal bundle $N_\F$ of $\F$ in $M$ splits into into its part in
$Y_+$, $N_\F':=N_\F\cap TY_+$, and the symplectic orthogonal of
$TY_+|\F$, which is canonically isomorphic to the trivial bundle
$\G/\T$. The weights for the $T$-action on $\G/\T$ are of course
simply the positive roots $\beta\in \T^*$ of $G$. Therefore, by
taking the trivial connection on $\G/\T$,
\[ \det(I- e^{-\f{i}{2\pi}F_\T(N_\F,\xi)})=\prod_{\beta>0}
(1-e^{-2\pi i \l \beta,\xi\r})\det(I-e^{-\f{i}{2\pi}
F_\T(N'_\F,\xi)}), \]
hence
\[ \chi(e^\xi)=\sum_{w\in W}\f{1}{\prod_{\beta>0}
(1-e^{-2\pi i \l \beta,w^{-1}(\xi)\r})}\sum_{\F\subset Y_+}
\int_\F \frac{{Td}\,(\F)e^{\omega
+ 2\pi i \l J_\F,w^{-1}(\xi)\r}}
{\det(I-e^{-\f{i}{2\pi}F_\T(N'_\F,w^{-1}(\xi))})}.\]
To the sum
\[ \sum_{\F\subset Y_+}\int_\F \frac{{Td}\,(\F)e^{\omega
+ 2\pi i \l J_\F,\xi\r}}
{\det(I-e^{-\f{i}{2\pi}F_\T(N'_\F,\xi)})},\]

we can apply the Localization Formula, this time in the opposite
direction, and get that it is equal to the above expression
$\chi'(e^\xi)$. This gives
\[\chi(e^\xi)=\sum_{w\in W} \f{\chi'(e^{w^{-1}(\xi)})}{  \prod_{\beta>0}
(1-e^{-2\pi i \l \beta,w^{-1}(\xi)\r})}.\]
But
\[   \prod_{\beta>0}
(1-e^{-2\pi i \l \beta,w^{-1}(\xi)\r})=\det(w)\,
e^{-2\pi i\l w(\delta)-\delta,\xi\r}\,\prod_{\beta>0}
(1-e^{-2\pi i \l \beta,\xi\r}),\]
where $\delta=\f{1}{2}\sum_{\beta>0}\beta$ is the magic weight.
Weyl's Character Formula hence shows that
\[ \chi(e^\xi)=\sum_{\mu\in\Lambda_+}
N'(\mu)\sum_{w\in W}
\det(w)\,\f{e^{2\pi i\l w(\delta+\mu)-\delta,\xi\r}}
{\prod_{\beta>0}
(1-e^{-2\pi i \l \beta,\xi\r})}=\sum_{\mu\in\Lambda_+}
N'(\mu) \,\chi_\mu(e^\xi),\]
where $\chi_\mu$ is the character of the irreducible representation
corresponding to $\mu$. This proves $N(\mu)=N'(\mu)$.
\bigskip\ebew

\noindent{\bf Remarks.}
\begin{enumerate}
\item
If $J(M)\not\subset \G^*_{reg}$, it is still possible to derive
a formula for $N(\mu)$ similar to (\ref{latticepoints}),
following part II of Guillemin-Prato \cite{GP90}.
However, this formula involves an additional ``shift'', so that
(\ref{11}) gets replaced by an equation of the form
\[ Al=m\nu+\sigma,\,l_j\ge l_j^0.\]
In general,
the number of integer solutions of such an equation is not an
arithmetic polynomial for all $m\in\N$, even though this is
true for large $m$.

\item On the other hand, Theorem \ref{stepw} does
not require that $\mu$ is a regular value of $J$. Even in the
singular case, it is therefore sufficient to prove Multiplicity
Formulas under the assumption $m>>0$.
\end{enumerate}

\noindent{\bf Acknowledgements.}

I would like to thank V. Guillemin, J. Kalkman, E. Lerman,
R. Sjamaar and C. Woodward for useful comments and discussions.
I am very much indebted to Victor Guillemin,
whose recent results \cite{G94} on Multiplicity Formulas of
Riemann-Roch type  have been the basic motivation for this work.
The stationary phase
version (\ref{jklocf}) of the Jeffrey-Kirwan Localization Theorem
was worked out jointly with Jaap Kalkman, who has also been
a great help in explaining equivariant cohomology to me.
This work was carried out when I was visiting scholar at the
M.I.T., and I wish to thank the Mathematics Department for its
hospitality.

\bigskip

\noindent{\bf Postscript:} After completing this article, we
learned about independent work of M. Vergne \cite{V94},
who has made a  different application of equivariant cohomology
to the multiplicity problem.

\end{document}